\documentclass{amsart}
\usepackage{amssymb}

\newtheorem{theorem}{Theorem}

\newtheorem{corollary}[theorem]{Corollary}

\newtheorem{remark}[theorem]{Remark}

\numberwithin{equation}{section}
\input{tcilatex}

\begin{document}
\title[Positive mass theorem]{Positive mass theorem and the Yamabe equation
on CR manifolds}
\author{Jih-Hsin Cheng}
\dedicatory{Dedicated to Professor Josip Globevnik on his 80th birthday}
\address{Institute of Mathematics, Academia Sinica and NCTS, Taipei Office,
Taipei, Taiwan, 10617, R.O.C.}
\email{cheng@math.sinica.edu.tw}
\urladdr{http://www.math.sinica.edu.tw}
\thanks{}

\begin{abstract}
Our goal is to survey the development of positive mass theorem and the
Yamabe equation on $CR$ manifolds in recent years. We introduce the notion
of the mass in several complex variables or $CR$ geometry. We then consider
the Yamabe problem on $CR$ manifolds to find a minimizer for the $CR$%
-Sobolev quotient. The positive mass theorem plays a key role in finding a
solution to the Yamabe equation with minimum energy for the positive
curvature case. We mainly focus on the team works in the following three
papers \cite{CMY17}, \cite{CMY23} and \cite{CC22}, on a positive mass
theorem in 3-dimensional $CR$ geometry, the $CR$-Sobolev quotient of Rossi
spheres, and the $5$-dimensional situation, respectively.
\end{abstract}

\maketitle

\section{Introduction}

We study a strongly pseudoconvex domain in several complex variables through
the geometry of its boundary. A higher-dimensional Riemann mapping theorem
(or when a strongly pseudoconvex domain is biholomorphic to the unit ball)
has been asked and there is obstruction on the geometry of the boundary-the
so called $CR$ (Cauchy-Riemann) geometry. In the old literature, a $CR$
structure was called a pseudoconformal structure. As a consequence of our
approach to the $CR$ Yamabe (minimizer) problem, we provide a single number
for such bounded domains of positive $CR$ Yamabe/Tanaka-Webster class, whose
vanishing characterizes the domain to be biholomorphic to the unit ball.
This number is called $p$-mass ("$p"$ stands for "pseudohermitian"). See (%
\ref{1.2}) and Corollary \ref{C-3-2} in Section \ref{Sec3}. For background
material, we refer the reader to Section \ref{Sec2}. The main results
discussed in this survey paper include two joint papers with Andrea
Malchiodi and Paul Yang: \cite{CMY17} on a positive mass theorem in
3-dimensional $CR$ geometry, \cite{CMY23} on the $p$-mass and the $CR$%
-Sobolev quotient of Rossi spheres, and another joint paper with Hung-Lin
Chiu \cite{CC22} on the case of 5-dimensional $CR$ spin manifolds.

On an odd dimensional manifold $M$, a contact structure $\xi $ is a natural
geometric structure to consider. Moreover, a contact manifold $(M,\xi )$
arising as the boundary of a domain in $\mathbb{C}^{n+1}$ (or a complex
manifold) inherits a complex structure, called $CR$ (Cauchy-Riemann)
structure. The $CR$ structure essentially reflects or controls the complex
structure of the inside domain. We can talk about abstract $CR$ structures
on a contact manifold $M$ (see Section \ref{Sec2}). We consider the
following CR Yamabe equation with critical Sobolev exponent (see (\ref{YE})
and notations in Section \ref{Sec2}):%
\begin{equation}
-(2+\frac{2}{n})\Delta _{b}u+Wu=u^{1+\frac{2}{n}}\text{ on }M.  \label{1-1}
\end{equation}%
\noindent Here $-\Delta _{b}$ and $W$ denote the (\textquotedblleft
positive") sublaplacian and the Tanaka-Webster scalar curvature,
respectively (see (\ref{SubL}) and (\ref{Wsc}) in Section \ref{Sec2}). There
is a variational formulation for equation (\ref{1-1}). Namely the energy is
provided by the following CR-Sobolev quotient (see (\ref{YMJ}) in Section %
\ref{Sec2}): 
\begin{equation}
Q(v):=\frac{\int_{M}((2+\frac{2}{n})|\nabla _{b}v|^{2}+Wv^{2})dV_{\theta }}{%
(\int_{M}v^{2+2/n}dV_{\theta })^{n/(n+1)}}\text{ for }v>0\text{ smooth.}
\label{1-2}
\end{equation}%
\noindent The main goal of this paper is to find a solution $u$ to (\ref{1-1}%
) with minimum energy $Q(u)$ equal to%
\begin{equation}
\inf_{v>0,\text{ smooth}}Q(v)=:\mathcal{Y}(M,J)  \label{1-3}
\end{equation}%
\noindent on a closed (compact with no boundary) $CR$ manifold provided $%
\mathcal{Y}(M,J)$ $>$ $0$. Let us give a brief history about this problem
below.

There has been a far-reaching analogy between conformal and $CR$ geometries.
Following the approaches used in conformal geometry by H. Yamabe, N.
Trudinger and T. Aubin, in 1987 D. Jerison and J. Lee \cite{JL87} showed the
analogous results in $CR$ geometry. That is, the $CR$ Yamabe constant $%
\mathcal{Y}(M,J)$ depends only on the $CR$ structure $J$ of $M$ and $%
\mathcal{Y}(M,J)$ $\leq $ $\mathcal{Y}(S^{2n+1},\hat{J})$, where $(S^{2n+1},$
$\hat{J})$ is the standard CR sphere with the induced CR structure $\hat{J}$
from $\mathbb{C}^{n+1}$. In addition, if $\mathcal{Y}(M,J)<\mathcal{Y}%
(S^{2n+1},\hat{J})$, then $\mathcal{Y}(M,J)$ is attained for some positive $%
C^{\infty }$ function $u$ (by the compactness of solutions to a family of
approximate equations), hence the $CR$ Yamabe minimizer problem for $(M,J)$
is solvable.

Recall that a $CR$ manifold is called spherical if it is locally $CR$
equivalent to the standard $CR$ sphere $(S^{2n+1},$ $\hat{J})$. In the case
that $2n+1\geq 5$ and ($M,J)$ is not spherical, Jerison and Lee \cite{JL89}
showed that $\mathcal{Y}(M,J)<\mathcal{Y}(S^{2n+1},\hat{J})$ by a test
function estimate. For the remaining cases: either (i) $\dim M=3$ or (ii) $%
\dim M$ $\geq $ $5$ and $M$ is spherical, we need a positive mass theorem to
show that $\mathcal{Y}(M,J)<\mathcal{Y}(S^{2n+1},\hat{J})$ unless $M$ is $CR$
equivalent to the standard $CR$ sphere. When $\dim M=3$, this was shown by
Malchiodi, Yang and myself \cite{CMY17}. See Theorem \ref{T-CMY17} in
Section \ref{Sec3}; the condition that the CR Paneitz operator of $M$ is
nonnegative turns out to be equivalent to the embeddability of $J$ \cite%
{Tak20}, see Remark \ref{R-CMY17}. As a consequence of Theorem \ref{T-CMY17}%
, we can solve the $CR$ Yamabe minimizer problem. See Theorem \ref{T-CMY17-2}%
.

There are non-standard $CR$ spheres, called Rossi spheres. In Section \ref%
{Sec4}, we show that their ($p$-)mass is negative, see Theorem \ref{T-CMY23}%
. Moreover, the infimum of the $CR$-Sobolev quotient of each Rossi sphere
coincides with that of the standard sphere and is not attained. See Theorem %
\ref{T-CMY23-2}.

Finally when $\dim M$ $\geq $ $5$ and $M$ is spherical, this was finished by
Chiu, Yang and myself (\cite{CCY14}; see Theorem \ref{T-CCY14} in Section %
\ref{Sec3}) through showing that the developing map is injective. However in
the case $\dim M$ $=$ $5$, we need an extra condition on the growth rate of
the Green's function on $\tilde{M},$ the universal cover of $M$. So in the
case $\dim M$ $=$ $5$, the $CR$ positive mass theorem is not really
completed. In the paper \cite{CC22}, we show that for $\dim M$ $=$ $5,$ $M$
being spherical, if in addition $M$ has a spin structure, then we have a $CR$
positive mass theorem, and hence the $CR$ Yamabe minimizer problem is
solvable (see Theorem \ref{YMP} in Section \ref{Sec5}).

For some comments on pseudohermitian Penrose inequality and other
interpretations of the $p$-mass, we refer the reader to \cite{Cheng17}.

\bigskip

{\large Acknowledgement.} The author would like to thank the Ministry of
Science and Technology of Taiwan for the support: grant no. MOST
114-2115-M-001-004. He would also like to thank Hung-Lin Chiu, Paul Yang and
Andrea Malchiodi for many discussions and constant interest in this
direction of research work. The material in this survey article was
presented in the workshop on Complex Analysis, Geometry, and Dynamics IV,
Portoro\v{z} 2025. The author would also like to thank the conference
organizers for the invitation to participate this meeting.

\section{Background material on $CR$ manifolds\label{Sec2}}

\subsection{Preliminary facts}

We introduce some basic notions and formulas in pseudohermitian geometry. We
refer the reader to N. Tanaka \cite{Ta75}, S. Webster \cite{Web78}, \cite%
{Lee86} and references therein.

Let $(M^{2n+1},\xi )$ denote a contact manifold with a coorientable (i.e. $%
TM/\xi $ is trivial) contact structure (or bundle) $\xi .$ A CR manifold $%
(M^{2n+1},\xi ,J)$ or $(M^{2n+1},J)$ (with $\xi $ suppressed) is a contact
manifold $(M^{2n+1},\xi )$ equipped with an almost complex structure, i.e.
an endomorphism $J:\xi \rightarrow \xi $ defined on $\xi $ such that $%
J^{2}=-1$. The endomorphism $J$ decomposes the complexification of $\xi $
into the direct sum of bundles of holomorphic vectors and anti-holomorphic
vectors $\xi \otimes \mathbb{C}=\xi _{1,0}\oplus \xi _{0,1}.$ We assume that 
$J$ is integrable, that is, $J$ satisfies the formal Frobenius condition $%
[\xi _{1,0},\xi _{1,0}]\subset \xi _{1,0}$ (as sections). A contact form $%
\theta $ is a global one-form such that $\xi =\ker {\theta }$ (exists by
coorientation of $\xi )$. A pseudohermitian manifold $(M^{2n+1},J,\theta )$
(with $\xi $ suppressed) is a contact manifold with a choice of CR structure 
$J$ together with a choice of contact form $\theta $. The Levi metric $%
L_{\theta }$ (assumed to be positive definite) is defined by 
\begin{equation}
L_{\theta }(X,Y):=\frac{1}{2}d\theta (X,JY)\text{ for all }X,Y\in \xi
\label{Levi-0}
\end{equation}%
(we use the convention that $\eta \wedge \vartheta (V,W)$ = $\eta
(V)\vartheta (W)$ $-$ $\eta (W)\vartheta (V)$ for 1-forms $\eta ,$ $%
\vartheta ,$ vectors $V,$ $W$). Let $T$ denote the Reeb vector field
associated to $\theta ,$ the unique vector field such that $\theta (T)$ $=$ $%
1$ and $L_{T}\theta $ $=$ $0$ ($L_{T}$ means the Lie derivative in the
direction $T).$ For a choice of (admissible) coframe $\theta ^{\alpha }$
with $\theta ^{\alpha }(T)$ $=$ $0,$ we have the Levi equation%
\begin{equation}
d\theta =ih_{\alpha \bar{\beta}}\theta ^{\alpha }\wedge \theta ^{\bar{\beta}%
}.  \label{Levi}
\end{equation}

In 1978, S. Webster \cite{Web78} (cf. an equivalent formulation in \cite%
{Ta75} by N. Tanaka) showed that there is a natural connection in the bundle 
$\xi _{1,0}$ adapted to a pseudohermitian structure $(J,\theta )$. Locally,
there exist unique $1$-forms $\theta _{\alpha }{}^{\beta }$ (connection
forms)$,\ \tau ^{\beta }$ (torsion forms) satisfying the structure equations%
\begin{eqnarray}
d\theta ^{\beta } &=&\theta ^{\alpha }\wedge \theta _{\alpha }\text{ }%
^{\beta }+\theta \wedge \tau ^{\beta },  \label{SE} \\
0 &=&\theta _{\alpha }\text{ }^{\beta }+\theta _{\bar{\beta}}\text{ }^{\bar{%
\alpha}},\text{ }0=\tau _{\beta }\wedge \theta ^{\beta }  \notag
\end{eqnarray}%
where \{$\theta ^{\beta }\}$ is a unitary coframe (meaning $h_{\alpha \bar{%
\beta}}=\delta _{\alpha \beta })$. Let $\{Z_{\beta }\}$ denote a unitary
frame of $\xi _{1,0}$ dual to \{$\theta ^{\beta }\}$. These forms $\theta
_{\alpha }{}^{\beta }$ satisfy the transformation law of connection forms,
so we can use them to define a connection. Let $T$ denote the Reeb vector
field associated to $\theta $. The \textbf{pseudohermitian }(or
Tanaka-Webster)\textbf{\ connection} $\nabla ^{p.h.}$ is defined by 
\begin{equation}
\begin{split}
\nabla ^{p.h.}{Z_{\alpha }}& =\theta _{\alpha }{}^{\beta }\otimes Z_{\beta }
\\
\nabla ^{p.h.}{Z_{\bar{\alpha}}}& =\theta _{\bar{\alpha}}{}^{\bar{\beta}%
}\otimes Z_{\bar{\beta}} \\
\nabla ^{p.h.}T& =0.
\end{split}
\label{Con_c}
\end{equation}

\noindent Differentiate the connection to define the curvature: $d\theta
_{\alpha }{}^{\beta }$ $-$ $\theta _{\alpha }{}^{\gamma }\wedge \theta
_{\gamma }{}^{\beta }$ $=$ $R_{\alpha }{}^{\beta }{}_{\rho \bar{\sigma}%
}\theta ^{\rho }\wedge \theta ^{\bar{\sigma}}$ $+\ $terms including the
torsion. The pseudohermitian-Ricci tensor is the hermitian form on $\xi
_{1,0}$ defined by 
\begin{equation*}
\rho (X,Y)=R_{\alpha \bar{\beta}}X^{\alpha }Y^{\bar{\beta}},
\end{equation*}%
where $X=X^{\alpha }Z_{\alpha },Y=Y^{\beta }Z_{\beta }$ and $R_{\alpha \bar{%
\beta}}=R_{\gamma }{}^{\gamma }{}_{\alpha \bar{\beta}}$. The Tanaka-Webster
scalar curvature is 
\begin{equation}
W:=R_{\beta }{}^{\beta },  \label{Wsc}
\end{equation}%
\noindent which is the contraction of the pseudohermitian-Ricci tensor. We
can also have a "real formulation" for the pseudohermitian structure $%
(J,\theta )$. Write $\theta ^{\alpha }$ $=$ $\omega ^{\alpha }+i\omega
^{n+\alpha }$ for real coframe fields $\{\omega ^{1},$ $..,$ $\omega ^{n},$ $%
\omega ^{n+1},$ $..,$ $\omega ^{2n}\}$ and $Z_{\alpha }$ $=$ $\frac{1}{2}%
(e_{\alpha }-ie_{n+\alpha })$ for real frame fields $\{e_{1},$ $..,$ $e_{n},$
$e_{n+1},$ $..,$ $e_{2n}\}$ (orthonormal with respect to the Levi metric $%
L_{\theta }).$ It is easily seen that $\{\omega ^{A}\}_{A=1,..,2n}$ is dual
to $\{e_{A}\}_{A=1,..,2n}.$ Write%
\begin{equation}
\nabla ^{p.h.}e_{A}=\omega _{A}\text{ }^{B}e_{B}  \label{Con_r}
\end{equation}

\noindent for real connection forms $\omega _{A}$ $^{B},$ $1$ $\leq $ $A,$ $%
B $ $\leq $ $2n$. Comparing (\ref{Con_c}) with (\ref{Con_r}) gives%
\begin{eqnarray}
\theta _{\alpha }{}^{\beta } &=&\omega _{\alpha }\text{ }^{\beta }+i\omega
_{\alpha }\text{ }^{n+\beta }\text{ and}  \label{Sym_r} \\
\omega _{\alpha }\text{ }^{n+\beta } &=&-\omega _{n+\alpha }\text{ }^{\beta }%
\text{, }\omega _{\alpha }\text{ }^{\beta }=\omega _{n+\alpha }\text{ }%
^{n+\beta }.  \notag
\end{eqnarray}

\noindent From the condition $0$ $=$ $\theta _{\alpha }$ $^{\beta }+\theta _{%
\bar{\beta}}$ $^{\bar{\alpha}}$ in (\ref{SE}) and (\ref{Sym_r}), it follows
that%
\begin{equation*}
\omega _{A}\text{ }^{B}+\omega _{B}\text{ }^{A}=0,\text{ }1\leq A,B\leq 2n.
\end{equation*}

\noindent Note that if we denote the scalar curvature associated to $\omega
_{A}$ $^{B}$ by $R,$ then we have%
\begin{equation}
W=\frac{1}{4}R.  \label{SC}
\end{equation}

\noindent Let $u_{\alpha \beta }$ denote the second covariant derivative of
a function $u$ in the directions $Z_{\alpha },Z_{\beta }.$ Define the
subgradient $\nabla _{b}$ and the sublaplacian $\Delta _{b}$ (or $\nabla
_{b}^{\theta }$ and $\Delta _{b}^{\theta }$ to indicate the dependence on $%
\theta )$ by%
\begin{equation*}
\nabla _{b}u:=u^{\alpha }Z_{\alpha }+u^{\bar{\alpha}}Z_{\bar{\alpha}},
\end{equation*}%
\begin{equation}
\Delta _{b}u:=u_{\alpha }\text{ }^{\alpha }+u_{\bar{\alpha}}\text{ }^{\bar{%
\alpha}}  \label{SubL}
\end{equation}%
\noindent (notice the \textquotedblleft positive" sign) where $u_{\alpha }$ $%
^{\alpha }$ :$=$ $u_{\alpha \bar{\beta}}h^{\alpha \bar{\beta}}$ $=$ $%
u_{\alpha \bar{\alpha}}$ for a unitary frame ($h^{\alpha \bar{\beta}}$ $=$ ($%
h_{\alpha \bar{\beta}})^{-1}$ $=$ $\delta _{\alpha \beta })$. Define the $CR$
invariant sublaplacian $L_{b}$ by%
\begin{equation}
L_{b}:=-b_{n}\Delta _{b}+W,\text{ }b_{n}=2+\frac{2}{n}.  \label{Lb}
\end{equation}%
\noindent Consider a new contact form $\hat{\theta}=u^{2/n}\theta $ for a
smooth positive function $u.$ $L_{b}$ rules the change of the Tanaka-Webster
scalar curvature:%
\begin{equation}
L_{b}u=\hat{W}u^{1+\frac{2}{n}}  \label{Lb-1}
\end{equation}%
\noindent where $\hat{W}$ is the Tanaka-Webster scalar curvature with
respect to $(J,\hat{\theta}).$ The Green's function $G_{p}$ of $L_{b}$ at $p$
satisfies%
\begin{equation}
L_{b}G_{p}=16\delta _{p}  \label{Lb-2}
\end{equation}%
\noindent where $\delta _{p}$ is the delta function with respect to the
volume form $dV_{\theta }$ $:=$ $\theta \wedge (d\theta )^{n}.$ We define
the $CR$ Yamabe constant $\mathcal{Y}(M,J)$ as follows: (cf. (\ref{1-3})) 
\begin{equation}
\mathcal{Y}(M,J):=\inf_{\hat{\theta}}\frac{\int_{M}\hat{W}dV_{\hat{\theta}}}{%
(\int_{M}dV_{\hat{\theta}})^{\frac{n}{n+1}}}=\inf_{0<u\in C^{\infty }(M)}%
\frac{\int_{M}(b_{n}|\nabla _{b}u|^{2}+Wu^{2})dV_{\theta }}{%
(\int_{M}u^{b_{n}}dV_{\theta })^{\frac{2}{b_{n}}}}  \label{YMJ}
\end{equation}

\noindent where $|\nabla _{b}u|^{2}:=2h^{\alpha \bar{\beta}}u_{\alpha }u_{%
\bar{\beta}}.$ Given a background $W$ with respect to $(J,\theta ),$ we aim
to find a solution $u$ to (\ref{Lb-1}) with $\hat{W}$ constant, say $1$ \
This is the so called Yamabe problem. The $CR$ Yamabe equation (with
critical Sobolev exponent) for $\hat{W}=1$ reads as follows (cf. (\ref{1-1}%
)): 
\begin{equation}
-b_{n}\Delta _{b}u+Wu=u^{1+\frac{2}{n}}.  \label{YE}
\end{equation}

The structure equations imply ($h_{\alpha \bar{\beta}}$ $=$ $\delta _{\alpha
\beta }$ for a unitary (co)frame)%
\begin{equation}
\begin{split}
\lbrack Z_{\bar{\beta}},Z_{\alpha }]& =ih_{\alpha \bar{\beta}}T+\theta
_{\alpha }{}^{\gamma }(Z_{\bar{\beta}})Z_{\gamma }-\theta _{\bar{\beta}}{}^{%
\bar{\gamma}}(Z_{\alpha })Z_{\bar{\gamma}}, \\
\lbrack Z_{\beta },Z_{\alpha }]& =\theta _{\alpha }{}^{\gamma }(Z_{\beta
})Z_{\gamma }-\theta _{\beta }{}^{\gamma }(Z_{\alpha })Z_{\gamma }, \\
\lbrack Z_{\alpha },T]& =A^{\bar{\gamma}}{}_{\alpha }Z_{\bar{\gamma}}-\theta
_{\alpha }{}^{\gamma }(T)Z_{\gamma },
\end{split}
\label{ide1}
\end{equation}

\noindent where we have written the torsion (forms) $\tau ^{\beta }$ $=$ $%
A^{\beta }$ $_{\bar{\alpha}}\theta ^{\bar{\alpha}}$ and $A^{\bar{\gamma}}$ $%
_{\alpha }$ $=$ $\overline{A^{\gamma }\text{ }_{\bar{\alpha}}}.$ Let $L_{T}$
denote the Lie differentiation in the direction $T.$ From (\ref{SE}) and the
third equality in (\ref{ide1}), it follows that%
\begin{equation}
L_{T}J=2iA^{\beta }\text{ }_{\bar{\alpha}}\theta ^{\bar{\alpha}}\otimes
Z_{\beta }-2iA^{\bar{\beta}}\text{ }_{\alpha }\theta ^{\alpha }\otimes Z_{%
\bar{\beta}}.  \label{LTJ}
\end{equation}

As a flat pseudohermitian manifold, the Heisenberg group plays an important
role in pseudohermitian geometry. We refer the reader to \cite{DT06}, \cite%
{Lee86}, \cite{Lee88} and \cite{Web78} for more details about
pseudohermitian geometry. Denote by $H_{n}$ the Heisenberg group, which is
the space $\mathbb{R}^{2n+1}$ with coordinates $(x_{\beta },y_{\beta },t)$
as a set. It is a $(2n+1)$-dimensional Lie group with group structure
defined by 
\begin{equation*}
(x,y,t)\circ (x^{\prime },y^{\prime },t^{\prime })=(x+x^{\prime
},y+y^{\prime },t+t^{\prime }+2yx^{\prime }-2xy^{\prime }).
\end{equation*}%
The associated Lie algebra is spanned by the following left invariant vector
fields 
\begin{equation}
\mathring{e}_{\beta }:=\frac{1}{\sqrt{2}}\left( \frac{\partial }{\partial
x^{\beta }}+2y^{\beta }\frac{\partial }{\partial t}\right) ,\ \ \mathring{e}%
_{n+\beta }:=\frac{1}{\sqrt{2}}\left( \frac{\partial }{\partial y^{\beta }}%
-2x^{\beta }\frac{\partial }{\partial t}\right) ,\ \ \mathring{T}:=\frac{%
\partial }{\partial t}.  \label{7-8a}
\end{equation}%
The associated standard CR structure $\mathring{J}$ and contact form $%
\mathring{\theta}$ (or denoted by $\Theta $) are defined respectively by%
\begin{eqnarray}
&\mathring{J}\mathring{e}_{\beta }{:=}\mathring{e}_{n+\beta },\text{ }%
\mathring{J}\mathring{e}_{n+\beta }:=-\mathring{e}_{\beta },&  \label{A1} \\
&\mathring{\theta}\text{ (or }\Theta \text{)}{:=} dt+\sum_{\beta
=1}^{n}(iz^{\beta }dz^{\bar{\beta}}-iz^{\bar{\beta}}dz^{\beta }).&  \notag
\end{eqnarray}%
Here $z^{\beta }$ $:=$ $x^{\beta }+iy^{\beta }$. The contact bundle is $%
\mathring{\xi}$ $:=\text{ ker }\mathring{\theta}$. We linearly extend $%
\mathring{J}:\mathring{\xi}\otimes \mathbb{C}\rightarrow \mathring{\xi}%
\otimes \mathbb{C}$. Let $\mathring{Z}_{\beta }:=\frac{1}{2}(\mathring{e}%
_{\beta }-i\mathring{e}_{n+\beta })=\frac{1}{\sqrt{2}}\left( \frac{\partial 
}{\partial z^{\beta }}+iz^{\bar{\beta}}\frac{\partial }{\partial t}\right) .$
Then for all $\beta ,\gamma ,$ $\mathring{J}\mathring{Z}_{\beta }=i\mathring{%
Z}_{\beta },\ \mathring{J}\mathring{Z}_{\bar{\beta}}=-i\mathring{Z}_{\bar{%
\beta}}$ and $[\mathring{Z}_{\beta },\mathring{Z}_{\gamma }]=0\ \left(
\Rightarrow \lbrack \mathring{\xi}_{1,0},\mathring{\xi}_{1,0}]\subset 
\mathring{\xi}_{1,0}\right) $ where $\mathring{\xi}\otimes \mathbb{C}=%
\mathring{\xi}_{1,0}\oplus \mathring{\xi}_{0,1}.$ It is easily seen that the
frame $\{\mathring{T},\mathring{Z}_{\beta },\mathring{Z}_{\bar{\beta}}\}$ is
dual to the coframe $\{\mathring{\theta},\sqrt{2}dz^{\beta },\sqrt{2}dz^{%
\bar{\beta}}\}$. If we regard $\{\mathring{e}_{\beta },\mathring{e}_{n+\beta
}|1\leq \beta \leq n\}$ as an orthonormal basis, then this defines a metric
on $\mathring{\xi}$, which equals the Levi metric $L_{\mathring{\theta}}$
given by $L_{\mathring{\theta}}(X,Y)=\frac{1}{2}d\mathring{\theta}(X,%
\mathring{J}Y)$ for all $X,Y\in \mathring{\xi}$. The standard
pseudohermitian connection on $H_{n}$ is defined by 
\begin{equation*}
\mathring{\nabla}^{p.h.}\mathring{e}_{\beta }=\mathring{\nabla}^{p.h.}%
\mathring{e}_{n+\beta }=\mathring{\nabla}^{p.h.}\mathring{T}=0.
\end{equation*}%
It follows that the pseudohermitian connection forms $\mathring{\theta}%
_{\alpha }$ $^{\gamma }$ vanish:%
\begin{equation}
\mathring{\theta}_{\alpha }\text{ }^{\gamma }=0  \label{A1-1}
\end{equation}%
We define the Heisenberg norm $\rho $ on $H_{n}$ by%
\begin{equation}
\rho ^{4}=(|z|^{4}+t^{2})  \label{A2}
\end{equation}%
where $|z|^{2}=\sum_{\beta =1}^{n}|z^{\beta }|^{2}$.

\subsection{Asymptotically Heisenberg manifolds and the $p$-mass}

A noncompact pseudohermitian ($2n+1)$-manifold $(N,J,\theta )$ ($J$ $:$ $CR$
structure, $\theta $ $:$ contact form) is asymptotically flat (or
asymptotically Heisenberg) if $N=N_{0}\cup N_{\infty },$ $N_{0}$ compact (or
complete), and there is a diffeomorphism: $N_{\infty }$ $\rightarrow $ ($%
\mathbb{C}^{n}\times \mathbb{R})\backslash B_{r},$ $r$ $>>$ $0,$ s.t.%
\begin{eqnarray}
\theta &=&(1+c_{n}A\rho ^{-2n}+O(\rho ^{-2n-1}))\Theta  \label{1.1} \\
&&+O(\rho ^{-2n-1})_{\beta }dz^{\beta }+O(\rho ^{-2n-1})_{\bar{\beta}}dz^{%
\bar{\beta}},  \notag
\end{eqnarray}

\begin{eqnarray*}
\theta ^{\alpha } &=&O(\rho ^{-2n-1})^{\alpha }\Theta +O(\rho ^{-2n-2})_{%
\bar{\beta}}^{\alpha }dz^{\bar{\beta}} \\
&&+[(1+c_{n}A\rho ^{-2n})\delta _{\beta }^{\alpha }+O(\rho ^{-2n-1})_{\beta
}^{\alpha }]\sqrt{2}dz^{\beta }
\end{eqnarray*}%
\noindent for some unitary coframe $\theta ^{\alpha },$ where $(z^{1},$ $%
z^{2},$ $...,$ $z^{n},$ $t$) $\in $ $C^{n}\times R$ (identified with
Heisenberg group $H_{n}$), $\Theta :=dt+iz^{\alpha }d\bar{z}^{\alpha }-i\bar{%
z}^{\alpha }dz^{\alpha }$ (cf. (\ref{A1}))$,$ $\rho $ $:=$ ($%
t^{2}+|z|^{4})^{1/4},$ and $B_{r}$ $:=$ $\{\rho $ $<$ $r\}$. The notation $%
O(\rho ^{-2n-1})^{\alpha }$ ($O(\rho ^{-2n-1})_{\beta }^{\alpha },$ resp.)
means the components are all $O(\rho ^{-2n-1}).$ We may also require the
Tanaka-Webster scalar curvature $R_{J,\theta }$ $\in $ $L^{1}(N,\theta
\wedge (d\theta )^{n}).$ This class of asymptotically Heisenberg manifolds
includes the blow-ups of closed spherical pseudohermitian manifolds (see
later comments). We can then define the $ADM$-like mass, called $p$-mass,
for such $(N,J,\theta )$ by%
\begin{equation}
m(J,\theta ):=\lim_{r\rightarrow \infty }ni\doint\limits_{\partial
B_{r}}\sum_{\gamma =1}^{n}\theta _{\gamma }{}^{\gamma }\wedge \theta \wedge
(d\theta )^{n-1}  \label{1.2}
\end{equation}

\noindent where $\theta _{\alpha }{}^{\beta }$'$s$ denote the connection
forms of pseudohermitian structure $(J,\theta )$ (cf. \cite{CMY17} for the $%
n $ $=$ $1$ case; \cite[Subsection 3.1]{CC22} for higher dimensional
situation)$.$ The $p$-mass is defined to kill the boundary term when we
compute the first variation of the pseudohermitian Einstein-Hilbert action:%
\begin{eqnarray*}
&&\delta _{J}\{-\int_{N}W_{J,\theta }\theta \wedge (d\theta )^{n}+m(J,\theta
)\} \\
&=&\frac{n}{2}\int_{N}<A_{J,\theta },\delta J>\theta \wedge (d\theta )^{n}
\end{eqnarray*}

\noindent where $W_{J,\theta }:=W$ is the Tanaka-Webster scalar curvature
(see (\ref{Wsc})) and $A_{J,\theta }$ is the torsion tensor (see more
details in \cite{CMY17} and (\ref{HE}) for the $n$ $=$ $1$ case):%
\begin{equation*}
A_{J,\theta }:=A^{\bar{\beta}}\text{ }_{\alpha }\theta ^{\alpha }\otimes Z_{%
\bar{\beta}}+A^{\beta }\text{ }_{\bar{\alpha}}\theta ^{\bar{\alpha}}\otimes
Z_{\beta }
\end{equation*}

\noindent (cf. (\ref{LTJ}); see \cite[(3.7) with $F_{\gamma }^{l}=0$]{ACMY24}%
).

Our goal is to show the nonnegativity of $m(J,\theta )$ under the condition $%
W_{J,\theta }$ $\geq $ $0$ and some other constraint(s). Also we want to
characterize the equality case $m(J,\theta )$ $=$ $0,$ namely, to show that $%
m(J,\theta )$ $=$ $0$ if and only if $(N,J,\theta )$ is isomorphic as
pseudohermitian manifold to the Heisenberg group $H_{n}$ with standard flat
pseudohermitian structure.

\subsection{Contact Dirac operator}

Let us start with a general spin$^{c}$ structure on a contact bundle $\xi $
over an asymptotically Heisenberg manifold $N$ of dimension $2n+1.$ Let $%
\mathcal{S}$ denote the spinor bundle with a spin$^{c}$ connection $\nabla $
compatible with the pseudohermitian connection $\nabla ^{p.h.}.$ Let $e_{1},$
$...,$ $e_{2n}$ denote an orthonormal basis of $\xi $ with respect to the
Levi metric. Denote the contact Dirac operator $D_{\xi }$ by%
\begin{equation}
D_{\xi }\psi =\sum_{\alpha =1}^{2n}\Gamma (e_{\alpha })\nabla _{e_{\alpha
}}\psi  \label{c0}
\end{equation}

\noindent (summation convention) for a section $\psi $ of $\mathcal{S},$
where $\Gamma $ denotes the Clifford multiplication. Let $T$ denote the Reeb
vector field associated to the contact form $\theta ,$ i.e., $\theta (T)$ $=$
$1$ and $T$ $\lrcorner $ $d\theta $ $=$ $0.$ Let $D_{\xi }^{\ast },$ $\nabla
^{\ast }$ denote the adjoint operator of $D_{\xi },$ $\nabla ,$
respectively. We then have the following formula 
\begin{eqnarray}
D_{\xi }^{\ast }D_{\xi }\psi &=&\sum_{\alpha =1}^{2n}\nabla _{e_{\alpha
}}^{\ast }\nabla _{e_{\alpha }}\psi -2\sum_{\alpha =1}^{n}\Gamma (e_{\alpha
}e_{n+\alpha })\nabla _{T}\psi  \label{c1} \\
&&+\sum_{a<b}\Gamma (e_{a})\Gamma (e_{b})R^{\nabla }(e_{a},e_{b})\psi  \notag
\end{eqnarray}

\noindent where $R^{\nabla }(e_{a},e_{b})$, the curvature operator, is
defined by%
\begin{equation*}
R^{\nabla }(e_{a},e_{b}):=\nabla _{e_{a}}\nabla _{e_{b}}-\nabla
_{e_{b}}\nabla _{e_{a}}-\nabla _{\lbrack e_{a},e_{b}]}.
\end{equation*}

\noindent Decompose $R_{ab}^{\nabla }$ := $R^{\nabla }(e_{a},e_{b})$ as a
sum of the trace free part $\mathring{R}_{ab}^{\nabla }$ and the trace part $%
2^{-n}tr_{\mathcal{S}}R_{ab}^{\nabla }.$ A standard deduction shows that%
\begin{equation*}
\sum_{a<b}\Gamma (e_{a})\Gamma (e_{b})\mathring{R}_{ab}^{\nabla }\psi =\frac{%
1}{4}R\psi
\end{equation*}

\noindent where $R$ denotes the (real version) Tanaka-Webster scalar
curvature (cf. (\ref{SC})). On the other hand, $2^{-n}tr_{\mathcal{S}%
}R_{ab}^{\nabla }$ $=$ $F_{A}(e_{a},e_{b})$ in which the curvature 2-form $%
F_{A}$ :$=$ $dA$ and $2A$ is the connection form of an associated line
bundle $L_{\Gamma }$ ($\det \mathcal{S}$ $=$ $L_{\Gamma }^{\otimes
2^{n-1}}). $ Therefore we can reduce (\ref{c1}) to 
\begin{eqnarray}
D_{\xi }^{\ast }D_{\xi }\psi &=&\sum_{\alpha =1}^{2n}\nabla _{e_{\alpha
}}^{\ast }\nabla _{e_{\alpha }}\psi -2\sum_{\alpha =1}^{n}\Gamma (e_{\alpha
}e_{n+\alpha })\nabla _{T}\psi  \label{c2} \\
&&+\frac{1}{4}R\psi +\rho (F_{A})\psi  \notag
\end{eqnarray}

\noindent where $\rho (F_{A})$ $=$ $\sum_{a<b}\Gamma (e_{a})\Gamma
(e_{b})F_{A}(e_{a},e_{b}).$

\section{A positive mass theorem in dimension $3\label{Sec3}$}

\subsection{Kohn's Laplacian and $CR$ Paneitz operator}

To deal with the $T$-derivative term in (\ref{c2}), we consider the
canonical spin$^{c}$ structure with $\mathcal{S}$ $=$ $\Lambda ^{0,\ast },$
the bundle of all $(0,q)$ forms. In particular, we take $\psi $ $=$ $\bar{%
\partial}_{b}u$ $=$ $u_{,\bar{\beta}}\theta ^{\bar{\beta}}$ (summation
convention throughout the remaining part)$,$ a $(0,1)$ form with components
being derivatives of a complex function $u.$ Then we have%
\begin{eqnarray*}
D_{\xi }\psi &=&(\bar{\partial}_{b}+\bar{\partial}_{b}^{\ast })\circ \bar{%
\partial}_{b}u \\
&=&\bar{\partial}_{b}^{\ast }\circ \bar{\partial}_{b}u=\square _{b}u
\end{eqnarray*}

\noindent where $\square _{b}$ $:=$ $\bar{\partial}_{b}\circ \bar{\partial}%
_{b}^{\ast }$ $+$ $\bar{\partial}_{b}^{\ast }\circ \bar{\partial}_{b}$ is
Kohn's Laplacian. Note that $\bar{\partial}_{b}^{\ast }u$ $=$ $0.$ To solve
the contact Dirac equation $D_{\xi }\psi $ $=$ $0$ for $\psi $ with a
suitably asymptotic behavior at infinity is reduced to solving $\square _{b}u
$ $=$ $0$ for $u$ with corresponding behavior at infinity. On the other
hand, we can compute the $T$-derivative term in (\ref{c2}) as follows:%
\begin{eqnarray*}
&&-2\sum_{\alpha =1}^{n}\Gamma (e_{\alpha }e_{n+\alpha })\nabla _{T}\psi  \\
&=&2i(n-2)\sum_{\beta =1}^{n}u_{,\bar{\beta}0}\theta ^{\bar{\beta}}
\end{eqnarray*}

\noindent Taking $L^{2}$ inner product with $\psi $ $=$ $\bar{\partial}_{b}u$
and making use of the equation $0$ $=$ $\square _{b}u$ $=$ $-2u_{,\bar{\gamma%
}\gamma },$ we get 
\begin{equation*}
-\frac{n-2}{2n}\int_{M}(u,Pu)\theta \wedge (d\theta )^{n}
\end{equation*}

\noindent modulo a boundary term, in which $P$ is the $CR$ Paneitz operator
defined by (all lower indices with respect to a unitary frame)%
\begin{equation}
Pu=4(u_{,\bar{\gamma}\gamma \beta }+niA_{\gamma \beta }u_{,\bar{\gamma}})_{,%
\bar{\beta}}  \label{3-a}
\end{equation}

Observe that for $n$ $\geq $ $2$ (dimension $2n+1$ $\geq $ $5)$, $P$ is
nonnegative (for closed $N$ and open $N$ with suitably decayed test
functions) (\cite{GL88}). For dimension $\geq 5,$ by assuming $N$ is
pseudo-Einstein, we can absorb the trace curvature term into the scalar
curvature term. So by further assuming $R$ $\geq $ $0,$ we can have the
nonnegativity of the $p$-mass (which we pick up from the boundary terms
after the integration of taking the inner product of (\ref{c2}) with $\psi $%
). The above argument relies on the existence of a solution to $\square
_{b}u $ $=$ $0$ with suitably asymptotic behavior near the infinity.

When $N$ arises as the blow-up of a closed pseudohermitian manifold $M$ by a
Green function $G_{p}$ of conformal sublaplacian at a point $p$, we can
apply the positive mass theorem for $N$ to solve the $CR$ Yamabe problem for 
$M.$ This situation occurs if $G_{p}$ has the following expansion near $p:$%
\begin{equation}
G_{p}=c_{n}\rho ^{-2n}+A_{p}+O(\rho )  \label{c3}
\end{equation}

\noindent where $\rho $ is the Heisenberg distance in $CR$ normal
coordinates. Note that \textquotedblleft $A_{p}"$ is a multiple of the $p$%
-mass defined for the blow-up $N.$ We observe that (\ref{c3}) holds for $n$ $%
=$ $1$ (dimension 3 case) and for $M$ being spherical of all dimensions.

For such spherical $M$ of dimension $\geq $ 5 (extra technical condition in
dimension 5) with positive $CR$ Yamabe/Tanaka-Webster class (i.e., the $CR$
Yamabe constant $\mathcal{Y}(M,J)$ $>$ $0,$ see (\ref{YMJ}) for the notation 
$\mathcal{Y}(M,J))$, we can prove a positive mass theorem for "$A_{p}";$ see
Theorem \ref{T-CCY14} below and hence find solutions of the $CR$ Yamabe
problem with minimal energy through another approach (\cite{CCY14}). We omit
describing this approach here, but only state the result as follows. \ 

Let $\tilde{M}$ denote the universal covering space of $M.$ Take $q$ $\in $ $%
\pi ^{-1}(p)$ where $\pi $ $:$ $\tilde{M}$ $\rightarrow $ $M$ is the natural
projection$.$ Let $\tilde{G}_{q}$ be a positive minimal Green's function
(see \cite{CCY14} for the definition) on $\tilde{M}$ with pole at $q$. Define%
\begin{equation*}
s(M):=\inf \{s:\int_{\tilde{M}\backslash U_{q}}\tilde{G}_{q}^{s}\theta
\wedge (d\theta )^{n}<\infty \}
\end{equation*}

\noindent where $U_{q}$ is a neighborhood of $q.$ The $CR$ invariant $s(M)$
measures the integrability of $\tilde{G}_{q}$ on $\tilde{M}.$

\begin{theorem}
\label{T-CCY14} ( \cite{CCY14})\textit{\ Let }$M$\textit{\ be a closed
spherical }$CR$\textit{\ manifold of dimension }$2n+1$\textit{\ with
positive }$CR$ Yamabe/\textit{Tanaka-Webster class. Then, for }$n\geq 3,$%
\textit{\ }$A_{p}$\textit{\ }$>$\textit{\ }$0$\textit{\ unless }$M$\textit{\
is the standard }$CR$\textit{\ sphere. In case }$n$\textit{\ }$=$\textit{\ }$%
2,$\textit{\ the same result also holds if we assume further }$s(M)$\textit{%
\ }$<$\textit{\ }$1.$
\end{theorem}

For the 5-dimensional spherical case, we can also follow the spinor method
to prove a positive mass theorem under a different, but geometric, condition
(for instance, $c_{1}(K)$ $=$ $0$ where $K$ is the canonical line bundle).
See Section \ref{Sec5}.

\subsection{A positive mass theorem in dimension 3}

In this subsection, we would like to discuss the case of dimension 3 in more
details. We refer for some results here to \cite{CMY17}.

We consider a closed (compact with no boundary) three dimensional
pseudohermitian manifold $(M,J,\theta )$ of \emph{positive CR
Yamabe/Tanaka-Webster class}$.$ This means that the first eigenvalue of the 
\emph{conformal sublaplacian} 
\begin{equation*}
L_{b}:=-4\Delta _{b}+W,
\end{equation*}%
\noindent is strictly positive. Here $\Delta _{b}$ stands for the
(\textquotedblleft positive") sublaplacian of $M$ and $W$ for the
Tanaka-Webster (scalar) curvature. The conformal sublaplacian has the
following covariance property under a conformal change of contact form 
\begin{equation*}
\hat{L}_{b}(\varphi )=u^{-\frac{Q+2}{Q-2}}L_{b}(u\varphi );\qquad \quad \hat{%
\theta}=u^{2}\theta ,
\end{equation*}%
\noindent where $Q=4$ is the \emph{homogeneous dimension} of the manifold.
The conformal sublaplacian rules the change of the Tanaka-Webster curvature
under the above conformal deformation through the following formula (cf. (%
\ref{Lb-1})) 
\begin{equation*}
-4\Delta _{b}u+Wu=\hat{W}u^{\frac{Q+2}{Q-2}},
\end{equation*}%
\noindent where $\hat{W}$ is the Tanaka-Webster curvature corresponding to
the pseudohermitian structure $(J,\hat{\theta})$. The positivity of the $CR$
Yamabe/Tanaka-Webster class is equivalent to the condition (cf. (\ref{YMJ}))%
\begin{equation}
\mathcal{Y}(J)(:=\mathcal{Y}(M,J)):=\inf_{\hat{\theta}}\frac{\int_{M}\hat{W}%
\text{ }\hat{\theta}\wedge d\hat{\theta}}{\left( \int_{M}\hat{\theta}\wedge d%
\hat{\theta}\right) ^{\frac{1}{2}}}>0,  \label{m.0}
\end{equation}%
\noindent where $\hat{\theta}$ is any contact form which annihilates the
underlying contact bundle $\xi $. Under the assumption $\mathcal{Y}(J)>0$ we
have that $L_{b}$ is invertible, so for any $p\in M$ there exists a Green's
function $G_{p}$ for which 
\begin{equation}
\left( -4\Delta _{b}+W\right) G_{p}=16\delta {_{p}}.  \label{3-2b}
\end{equation}
\noindent One can show that in $CR$ normal coordinates $(z,t)$ \cite{JL89}, $%
G_{p}$ admits the following expansion 
\begin{equation}
G_{p}=\frac{1}{2\pi }\rho ^{-2}+A+O(\rho )  \label{3-2a}
\end{equation}
\noindent where $A$ is some real constant and where we have set $\rho
^{4}(z,t)=|z|^{4}+t^{2}$, $z\in \mathbb{C},t\in \mathbb{R}$. Having in mind
the Riemannian construction for the blow-up of a closed manifold, we
consider the new pseudohermitian manifold with a blow-up of contact form%
\begin{equation}
N=(M\setminus \{p\},J,\theta =G_{p}^{2}\hat{\theta}),  \label{m.1}
\end{equation}%
\noindent where $\hat{\theta}$ is suitably chosen. With an \emph{inversion
of coordinates} we then obtain a pseudohermitian manifold which has
asymptotically the geometry of the Heisenberg group. Starting from this
model, we give a definition of asymptotically flat pseudohermitian manifold
(see (\ref{1.1}) for the $n$ $=$ $1$ case) and we introduce its \emph{%
pseudohermitian mass} ($p$-mass) by the formula (cf. (\ref{1.2}) for the $n$ 
$=$ $1$ case) 
\begin{equation*}
m(J,\theta ):=i\oint_{\infty }\theta _{1}{}^{1}\wedge \theta
:=\lim_{L\rightarrow +\infty }i\oint_{S_{L}}\theta _{1}{}^{1}\wedge \theta ,
\end{equation*}%
\noindent where we have set $S_{L}=\left\{ \rho =L\right\} $, $\rho
^{4}=|z|^{4}+t^{2}$, and where $\theta _{1}{}^{1}$ stands for the connection
form of the structure. The above quantity is indeed a natural candidate,
since it satisfies a property analogous to the situation in general
relativity:%
\begin{equation}
\frac{d}{ds}|_{s=0}(-\int_{N}W_{J(s),\theta }\theta \wedge d\theta
+m(J(s),\theta ))=\int_{N}(A_{11}E_{\bar{1}\bar{1}}+A_{\bar{1}\bar{1}%
}E_{11})\theta \wedge d\theta  \label{HE}
\end{equation}%
\noindent where $\dot{J}=2E=2E_{11}\theta ^{1}\otimes Z_{\bar{1}}+2E_{\bar{1}%
\bar{1}}\theta ^{\bar{1}}\otimes Z_{1}$ and $A_{11}$ ($=A^{\bar{1}}$ $_{1}$)
denotes the pseudohermitian torsion with respect to unitary frame $Z_{1}$
and coframe $\theta ^{1}$. Moreover it coincides with the zeroth order term
in the expansion of the Green's function (\ref{3-2a}) for $L_{b}$:%
\begin{equation}
m(J,\theta )=48\pi ^{2}A.  \label{3-4a}
\end{equation}

We prove an integral formula for the p-mass, in the spirit of \cite{Wit81}.
To state this formula we need to introduce another conformally covariant
operator, the $CR$ Paneitz operator (cf. (\ref{3-a}))%
\begin{equation}
P\varphi :=4(\varphi {_{,\bar{1}}}^{\bar{1}}{_{1}}+iA_{11}\varphi
,^{1}),^{1}.  \label{3-4b}
\end{equation}
\noindent The operator $P$ satisfies the covariance property%
\begin{equation}
P_{(J,\hat{\theta})}=u^{-4}P_{(J,\theta )};\qquad \quad \hat{\theta}%
=u^{2}\theta  \label{3-4c}
\end{equation}%
\noindent (see \cite{Hi93}, \cite{GL88}). We prove then the following
integral formula, which holds for an asymptotically flat pseudohermitian
manifold $N$%
\begin{equation}
\frac{2}{3}m(J,\theta )=-\int_{N}|\Box _{b}\beta |^{2}\theta \wedge d\theta
+2\int_{N}|\beta _{,\overline{1}\overline{1}}|^{2}\theta \wedge d\theta
+2\int_{N}W|\beta _{,\overline{1}}|^{2}\theta \wedge d\theta +\frac{1}{2}%
\int_{N}\overline{\beta }P\beta \,\theta \wedge d\theta .  \label{m.2}
\end{equation}%
\noindent Here $\beta :N\rightarrow \mathbb{C}$ is a function satisfying 
\begin{equation}
\beta =\overline{z}+\beta _{-1}+O(\rho ^{-2+\varepsilon })\quad \text{near }%
\infty ;\qquad \qquad \square _{b}\beta =O(\rho ^{-4}),  \label{m.3}
\end{equation}%
\noindent with $\square _{b}=-2\beta _{,\overline{1}1}$ and with $\beta
_{-1} $ a suitable function with homogeneity $-1$ in $\rho $.

In the following theorem we give some general conditions which ensure the
nonnegativity of the $p$-mass, characterizing also the zero case as ($CR$
equivalent to) the standard three dimensional $CR$ sphere.

\begin{theorem}
\label{T-CMY17} \textit{Let }$M$\textit{\ be a smooth, strictly pseudoconvex
compact (with no boundary) CR manifold of dimension 3. Suppose }$\mathcal{Y}%
(J)>0,$\textit{\ and that the CR Paneitz operator is nonnegative (which is
equivalent to }$J$ being embeddable; see Remark \ref{R-CMY17} below) \textit{%
. Let }$p$\textit{\ }$\in $\textit{\ }$M$\textit{\ and let }$\theta $\textit{%
\ be a blow-up of contact form as in (\ref{m.1}). Then}
\end{theorem}

\textit{(a) }$m(J,\theta )\geq 0;$

\textit{(b) if }$m(J,\theta )=0,$\textit{\ }$M$\textit{\ is CR equivalent to 
}$S^{3},$\textit{\ endowed with its standard CR structure.}

\begin{remark}
\label{R-CMY17} The assumptions we give in Theorem \ref{T-CMY17} are
conformally invariant, and are needed to ensure the positivity of the
right-hand side in (\ref{m.2}). By the result in \cite{ChaCY12}, the
conditions on $\mathcal{Y}(J)$ and $P$ imply the embeddability of $M$.
Therefore $\bar{\partial}_{b}$ has closed range. It is then used to solve $%
\square _{b}\beta $ $=$ $0$ with $\beta $ satisfying (\ref{m.3}) and to deal
with the Paneitz term when converted back to $M$ (\cite{HY15}). Conversely,
the embeddability of $M$\ in fact implies the nonnegativity of $P$. See Yuya
Takeuchi's elegant proof in \cite{Tak20}. So the embeddability is equivalent
to $P$ $\geq $ $0$ under the condition $\mathcal{Y}(J)$ $>$ $0$.
\end{remark}

As a consequence of Theorem \ref{T-CMY17} and Remark \ref{R-CMY17}, we have
the following results.

\begin{corollary}
\label{C-3-1} The $CR$ Yamabe equation (\ref{1-1}) (with $n=1$) has a
solution with minimum energy for $(M,J)$ embeddable.
\end{corollary}

\begin{corollary}
\label{C-3-2} (a version of generalized Riemann mapping theorem) Let $\Omega
\subset \mathbb{C}^{2}$ be a \textit{strictly pseudoconvex} domain close
enough to the unit ball $B^{2}\subset \mathbb{C}^{2}.$ Suppose $m(J,\theta
)=0.$ Then $\Omega $ is biholomorphic to $B^{2}.$
\end{corollary}

Our next main goal is to apply Theorem \ref{T-CMY17} to the study of the $CR$
Yamabe problem, namely finding conformal changes of contact form in order to
obtain constant Tanaka-Webster curvature. As for the classical Yamabe
problem, the cases $\mathcal{Y}(J)\leq 0$ are rather easy to deal with,
while the case $\mathcal{Y}(J)>0$ is the most difficult one. Calling $%
\mathcal{Y}_{0}$ the quotient for the standard $CR$ three sphere, by a
result in \cite{JL87} one always has%
\begin{equation*}
\mathcal{Y}(J)\leq \mathcal{Y}_{0},
\end{equation*}%
and if strict inequality holds then the problem is solvable. The strict
inequality is needed to ensure compactness of the minimizing sequences in (%
\ref{m.0}). This condition was verified in \cite{JL89} for (real) dimension
greater or equal to five, and for nonlocally spherical structures, in the
spirit of \cite{Aub76}. Our next result gives the strict inequality in the
three dimensional case, if $M$ is not $CR$ equivalent to the standard $CR$
3-sphere, under the same assumptions as in Theorem \textit{\ref{T-CMY17}}.

\begin{theorem}
\label{T-CMY17-2} \textit{Suppose we are under the assumptions of Theorem %
\ref{T-CMY17} \ Then either }$M$\textit{\ is the standard CR 3-sphere or if }%
$M$\textit{\ is not CR equivalent to the standard CR 3-sphere, one has }$%
\mathcal{Y}(J)<\mathcal{Y}_{0}.$\textit{\ In both cases, the }$CR$-Sobolev%
\textit{\ quotient (\ref{1-2}) admits a smooth minimizer.}
\end{theorem}

The $CR$ Yamabe problem for the case of three-dimensional $CR$ manifolds and
for spherical $CR$ manifolds was solved in \cite{Ga01} and \cite{GY01}
respectively. While the proof in these papers relies on topological
arguments, in the spirit of \cite{BaBr96}, our argument is based on direct
minimization and gives an extra variational characterization on the
solutions. To prove strict inequality we follow Schoen's argument in \cite%
{Sc84}, finding test functions which resemble a $CR$ bubble at a small
scale, and the Green's function $G_{p}$ at a larger one. More in general,
the analysis of the Yamabe problem in the $CR$ case has been so far less
precise than the Riemannian case: for example a basic difficulty is the lack
of a moving plane method, which is useful in general to derive a priori
estimates and to classify entire solutions.

\section{Rossi spheres, negative $p$-mass and $CR$-Sobolev quotient\label%
{Sec4}}

In this section, we would like to give a brief introduction to results in 
\cite{CMY23}.

In \cite[JL87]{JL87}, the counterpart of the result in \cite{Aub76} was
obtained, i.e. if the infimum of the CR-Sobolev quotient satisfies (cf. (\ref%
{1-3}))%
\begin{eqnarray*}
\mathcal{Y}(M,J) &:&=\inf_{\hat{\theta}}\frac{\int_{M}W_{\hat{\theta}}\,\hat{%
\theta}\wedge d\hat{\theta}}{\left( \int_{M}\hat{\theta}\wedge d\hat{\theta}%
\right) ^{\frac{1}{2}}} \\
&=&\inf_{u\in C^{\infty }(M),u>0}\frac{\int_{M}(4|\nabla _{b}u|^{2}+W_{{%
\theta }}u^{2})\,{\theta }\wedge d{\theta }}{\left( \int_{M}u^{4}{\theta }%
\wedge d{\theta }\right) ^{\frac{1}{2}}}<\mathcal{Y}(S^{3},J_{S^{3}}),
\end{eqnarray*}
\noindent then it is attained and a solution of (\ref{Lb-1}) exists (indeed,
this holds true in any dimension). The same authors verified this condition
when the dimension is greater or equal to five and $(M,J)$ is not \emph{%
spherical}, see \cite{JL89} and \cite{JL88}.

However, in the $CR$ setting new phenomena appear, related to the fact that
most three-dimensional structures are nonembeddable, differently from the
higher-dimensional case, see \cite{BdM75}, \cite{BE90}. In \cite{CMY17} some
results in the above directions were proved, assuming some global conditions
related to the embeddability of the abstract $CR$ structure.

As it happens in the Riemannian case, the ($p$-)mass is related to the
expansion of the Green's function of the conformal sublaplacian $L_{b}$ on a
compact manifold $M$. When $\mathcal{Y}(M,J)>0$ the latter operator is
invertible, so for any $p\in M$ there exists a Green's function $G_{p}$
verifying distributionally (cf. (\ref{3-2b}); here we use an extra factor $%
4\pi $ in the definition of $G_{p})$ 
\begin{equation*}
\left( -4\Delta _{b}+W\right) G_{p}=64\,\pi \,\delta _{p},
\end{equation*}%
where $\delta _{p}$ in the right-hand side stands for the Dirac delta w.r.t.
the volume measure $\theta \wedge d\theta $. In \emph{CR normal coordinates} 
$(z,t)$ (introduced in \cite{JL89} and discussed in \cite[Section 2]{CMY23}) 
$G_{p}$ writes as 
\begin{equation}
G_{p}=2\rho ^{-2}+A+O(\rho ),  \label{eq:Green-A}
\end{equation}%
for some $A\in \mathbb{R}$ and where $\rho ^{4}(z,t)$ is as above. For the
latter expansion, we refer to Proposition 5.2 in \cite{CMY17} and to \cite[%
Subsection 2.1]{CMY23} for our notation $O(\rho )$. Given $(M,J,\theta )$
compact and $p\in M$, consider a blow-up of contact form as follows 
\begin{equation}
N=(M\setminus \{p\},J,G_{p}^{2}\theta ).  \label{eq:bbuu}
\end{equation}%
As it is shown in \cite{CMY17}, via an inversion of coordinates, the
manifold $N$ turns out to have asymptotically the geometry of the Heisenberg
group, and its pseudohermitian mass satisfies (cf. (\ref{3-4a})) 
\begin{equation}
m=12\pi A  \label{eq:mass-A}
\end{equation}%
(see Lemma 2.5 in \cite{CMY17}\ there, and recall the difference of $4\pi $
in our current notation), where $A$ is as above. Using crucially a result in 
\cite{HY15}, in the same paper it was also proved that the pseudohermitian
mass is nonnegative (and zero only when $(M,J,\theta )$ is $CR$ equivalent
to $S^{3}$), provided that the \emph{CR Paneitz operator} $P$ (see (\ref%
{3-4b})) on $(M,J)$ is nonnegative definite. It has a relation to the $\log $%
-term coefficient in the Szeg\"{o} kernel expansion, and it is
pseudohermitian-covariant, namely $P_{\hat{\theta}}\varphi =e^{4f}P_{\theta
}\varphi $ (recall (\ref{3-4c})) for the conformal change $\theta =e^{2f}%
\hat{\theta}$. By a result in \cite{ChaCY12}, manifolds for which $P$ is
nonnegative and $R>0$ can be embedded into some $\mathbb{C}^{N}.$

The assumption on the positivity of the Paneitz operator is not technical,
as in \cite{CMY17} some counterexamples for the positivity of the
pseudohermitian mass were also given for structures (arbitrarily) close to
the spherical one, and hence with positive Tanaka-Webster scalar curvature.
In \cite{Tak20}, the positivity of the Paneitz operator is shown to hold for
embeddable $(M,J)$.

In \cite{CMY23}, we are concerned with \emph{Rossi spheres}: these are a
one-parameter-family of $CR$ structures on the 3-sphere of the form $%
S_{s}^{3}:=(S^{3},J_{(s)},\hat{\theta})$, where $\hat{\theta}$ is the
standard contact form on $S^{3}$ as in \cite[(1.4)]{CMY23}, and where $%
J_{(s)}$ is characterized by 
\begin{equation}
J_{(s)}Z_{1(s)}=iZ_{1(s)};\text{ \ }Z_{1(s)}=Z_{1}+\frac{s}{\sqrt{1+s^{2}}}%
Z_{\bar{1}},\text{ }Z_{\bar{1}(s)}=Z_{\bar{1}}+\frac{s}{\sqrt{1+s^{2}}}Z_{1}.
\label{eq:Js}
\end{equation}%
Rossi spheres are interesting because they are simple examples of $CR$
structures on the three-sphere that cannot be embedded in $\mathbb{C}^{N}$.
In \cite{Bur79}, it was shown that all the holomorphic functions on such
structures are even functions if $s\neq 0$. On the other hand, there are
explicit embeddings in $\mathbb{C}^{3}$ of the quotient of the Rossi spheres
by the antipodal map, see \cite{CS01}. By the above discussion, it follows
that the Paneitz operator cannot be nonnegative here. In addition, this
family of $CR$ structures are homogeneous (see \cite{CMY24} for more
information about Rossi spheres) and it is \emph{pseudo-Einstein}, i.e. $%
R_{,1}-iA_{11,\bar{1}}=0$, see \cite{CY13} as well as our notation for
covariant derivatives in \cite[Section 2.1]{CMY23}.

Our first main result in \cite{CMY23} is the following theorem.

\begin{theorem}
\label{T-CMY23} For $|s|$ small, s$\neq $0, the pseudohermitian mass of the
Rossi spheres $S_{s}^{3}$ is negative. More precisely, one has the expansion%
\begin{equation*}
m_{s}=-18\pi s^{2}+o(s^{2})\text{ \ for }s\simeq 0.
\end{equation*}
\end{theorem}

We remark that (a) we can generalize the construction of Rossi spheres in
Theorem \ref{T-CMY23} as follows. According to \cite[Proposition 3.3]{Fal92}%
, there exist deformations of the standard CR structure on $S^{3}/\Gamma $ ($%
\Gamma $=$\mathbb{Z}_{2}$ for the case of Rossi spheres), whose universal
covers are not embeddable. These CR structures (i.e., universal covers) are
likely to have negative mass.

(b) We can embed $S_{s}^{3}$/$\mathbb{Z}_{2}$ into $\mathbb{C}^{3}$ (see,
for instance, \cite{CS01}). So according to \cite{Tak20}, the $CR$-Paneitz
operator P on $S_{s}^{3}/\mathbb{Z}_{2}$ is nonnegative definite. On the
other hand, P on $S_{s}^{3}$ cannot be nonnegative definite by Theorem 1.1
and the positive mass theorem in \cite{CMY17} for $|s|$ small, $s\neq 0$, so
that the Tanaka-Webster scalar curvature of $S_{s}^{3}$ is positive. Thus,
for $|s|$ small, $s\neq 0$, $S_{s}^{3}/\mathbb{Z}_{2}$ provides an example
of $CR$ manifold having nonnegative definite P while its covering space $%
S_{s}^{3}$ does not have nonnegative definite P, answering a question raised
by Ngaiming Mok in a conference held in Hong Kong, 2014.

We saw before (in both low-dimensional Riemannian and $CR$ cases) that
positivity of the mass implies attainment of the Sobolev quotient. We also
strengthen the relation between mass and quotient by means of the following
result, which is in striking contrast with the Riemannian case.

\begin{theorem}
\label{T-CMY23-2} For $|s|$ small, s$\neq $0, the infimum of the CR-Sobolev
quotient of $S_{s}^{3}$ coincides with $\mathcal{Y}(S^{3},J_{S^{3}})$ and is
not attained.
\end{theorem}

We remark that (a) in \cite{Ga01} and \cite{GY01} the $CR$-Yamabe problem
was solved for every three dimensional $CR$ manifold, but their solutions
were found via variational arguments and they are not of minimal type.
Theorem \ref{T-CMY23-2} shows that the use of such methods may not be in
some cases necessary; (b) the phenomenon in Theorem \ref{T-CMY23-2} is
typical of some critical problems in a PDE context, like the Yamabe equation
on Euclidean domains with Dirichlet boundary conditions or the case of some
general elliptic operators on manifolds. However, to our knowledge this is
the first time this is displayed in a purely geometric smooth context.

Determining or estimating the \emph{mass} of a manifold is in general a hard
problem, since this is deeply related to the Green's function of the
conformal sublaplacian, which is a \emph{global} object. The mass also
appears as its zero-th order coefficient after a proper choice of conformal
representative and local coordinates. After recalling some preliminary facts
in \cite[Section 2]{CMY23} on $CR$ normal coordinates (introduced in \cite%
{JL89} and suited for the above expansion) and on Rossi spheres, we
specialize in \cite[Section 3]{CMY23} to the latter manifolds. For doing
this we need first to derive a suitable conformal factor satisfying a list
of conditions, and then express pseudo-hermitian coordinates depending on $s$%
. By the special expression of the Green's function in these coordinates, we
are able to determine it quite precisely near the north pole, up to the
constant term $A$ appearing in \eqref{eq:Green-A}. However, as we remarked
before, also some global features of the Green's function have to be
understood.

For doing this, by a Taylor expansion in $s$ worked-out at the beginning of 
\cite[Section 4]{CMY23} it is possible to characterize formally the Green's
function for the conformal sublaplacian on Rossi spheres up to an order $%
O(s^{3})$. One problem with this expansion is that it generates singular
terms, with a particularly bad behavior near the pole, if expressed with
respect to the standard complex coordinates of $\mathbb{C}^{2}$, where $%
S^{3} $ embeds. Also in this case non local terms appear, which we are able
to evaluate at the pole via some integral formula.

Via a careful analysis of all terms of order $1$, $s$ and $s^{2}$, we verify
then in the second part of the section that the global singular expansion on 
$S^{3}$ matches with the one done in $CR$ normal coordinates up to an order $%
O(s^{3})$. This allows us to prove Theorem \ref{T-CMY23}.

In \cite[Section 5]{CMY23}, arguing by contradiction, we analyse the
possible behaviours of minimizers for the $CR$-Sobolev quotient. Due to a
non-degeneracy result from \cite{MU02}, the analysis of minimizers can be
reduced to a finite-dimensional one, and we show that the $CR$-Sobolev
quotient of all candidate minimizers is strictly above the spherical one,
i.e. $\mathcal{Y}(S^{3},J_{S^{3}})$. With negative mass, this is expected
for highly concentrated profiles, reversing the expansion in \cite{Sc84}:
however such a property has to be obtained in all cases, i.e. even for
non-concentrated profiles, in order to guarantee that the infimum of the $CR$%
-Sobolev quotient is not attained. In \cite[Proposition 5.5]{CMY23} this is
proved for $s$ small in a fixed compact set of the $CR$ maps of $S^{3}$.
This is done starting with the expansion of the quotient on Rossi spheres
over the extremals of the quotient on the standard $S^{3}$, adding to them a
correction term that improves their accuracy as approximate critical points
for $s$ nonzero. One needs then to analyze the quotient in a regime with
loss of compactness, which is particularly delicate due to the following
reason. It is known from \cite{Sc84} that the mass of a (given) manifold
plays a role in the expansion for Sobolev quotients of highly concentrated
functions. In our case this must be done uniformly in $s$, and the problem
could be that the \emph{principal term} coming from the mass could become
negligible as $s\rightarrow 0$. To solve this issue we exploit a symmetry $%
s\rightarrow -s$ for Rossi spheres, discussed in \cite[Section 2]{CMY23},
which implies that all variational expansions are indeed even in $s$ and
hence the mass, which vanishes with $s$, gives still a dominant sign to the
asymptotic expansion of the $CR$-Sobolev quotient. Two appendices in \cite%
{CMY23} are devoted to the estimates of the latter quantity in two different
scaling regimes. To make the above arguments rigorous, we employ a
finite-dimensional reduction of the problem, via a fixed point argument,
which allows to solve for the $CR$-Yamabe equation on Rossi spheres up to a
Lagrange multiplier. We obtain in this way a manifold of approximate
solutions containing by construction all possible minimizers: our expansion
shown then that on this manifold the $CR$-Sobolev quotient is strictly
higher than $\mathcal{Y}(S^{3},J_{S^{3}})$, yielding our result.

\section{The Yamabe equation on $5$-dimensional contact spin manifolds\label%
{Sec5}}

In the paper \cite{CC22}, we show that for $\dim M$ $=$ $5,$ $M$ being
spherical, if in addition $M$ has a spin structure, then we have a $CR$
positive mass theorem, and hence the $CR$ Yamabe minimizer problem is
solvable (see Theorem \ref{YMP} below). For an asymptotically flat
pseudohermitian manifold $(N,J,\theta )$ (see (\ref{1.1})), we can talk
about the p-mass $m(J,\theta )$ (see (\ref{1.2})).

\begin{theorem}
\label{PMT} Suppose that $(N,J,\theta )$ is an asymptotically flat,
pseudohermitian and spin manifold of dimension 5. Assume that $J$ is
spherical and $(N,J,\theta )$ has the Tanaka-Webster scalar curvature $W\geq
0$. Then the p-mass $m(J,\theta )\geq 0$. Moreover, $m(J,\theta )=0$ if and
only if $(N,J,\theta )$ is isomorphic to the Heisenberg group $(H_{2},%
\mathring{J},\mathring{\theta})$.
\end{theorem}

\begin{corollary}
\label{PMT'} Suppose that $(M,\xi )$ is a closed (compact with no boundary),
contact and spin manifold of dimension 5. Assume that $J$ is a spherical $CR$
structure on $(M,\xi )$ with $\mathcal{Y(}M,J)$ $>$ $0.$ Then the associated
p-mass $m(J,\theta )\geq 0$. Moreover, $m(J,\theta )=0$ if and only if $%
(M,J) $ is $CR$ equivalent to the standard $CR$ $5$-sphere.
\end{corollary}

The proof of Theorem \ref{PMT} is based on a Weitzenbock-type formula:%
\begin{equation}
D_{\xi }^{2}=\nabla ^{\ast }\nabla +W-2\sum_{\beta =1}^{n}e_{\beta
}e_{n+\beta }\nabla _{T}  \label{1-4}
\end{equation}

\noindent where $D_{\xi }$ and $\nabla $\ denote the contact Dirac operator
and spin connection respectively\ (see (\ref{c0}) or \cite[(2.14)]{CC22} for
more details). The term involving $\nabla _{T}$ ($T$ is the Reeb vector
field associated to the contact form $\theta $) causes difficulty to solve
the Dirac equation $D_{\xi }^{2}\psi $ $=$ $0$ in general. However, in the
case of dimension 5 ($n=2$) we observe the following algebraic fact for
Clifford multiplication:%
\begin{equation}
\sum_{\beta =1}^{2}e_{\beta }e_{2+\beta }=0\text{ on }S^{+}(2n)\overset{n=2}{%
=}S^{+}(4)  \label{1-5}
\end{equation}

\noindent (not true on $S^{-}(4))$ where $S^{+}(2n)$ denotes the space of
positive spinors (see \cite[(2.3)]{CC22}). So for the dimension equal to 5,
the last term in (\ref{1-4}) disappears when acting on (sections of) $%
\mathbb{S}^{+}$ (the bundle of positive spinors) by (\ref{1-5}) as $\nabla
_{T}$ leaves $\mathbb{S}^{+}$ invariant$.$ It follows that $D_{\xi }^{2}$ is
subelliptic on $\mathbb{S}^{+}$ and hence we can find a spinor field $\psi $ 
$\in $ $\mathbb{S}^{+}$ such that $D_{\xi }^{2}\psi =0$ and $\psi $ tends to
a constant spinor at infinity (see \cite[Corollary 4.2]{CC22}). Applying (%
\ref{1-4}) to this spinor field $\psi $ and integrating after taking the
inner product with $\psi ,$ we then pick up (a positive multiple of) the
p-mass $m(J,\theta )$ from the boundary integral and obtain a Witten-type
formula for $m(J,\theta )$ (see \cite[(4.24) and (3.13)]{CC22}). So
nonnegativity of $m(J,\theta )$ follows. Note that we may not have $D_{\xi
}\psi $ $=$ $0$ which is usually used in deducing the Witten-type formula.
Instead we show that $D_{\xi }\psi $ has a fast enough decay order for our
purpose (see Proposition 4.3 of \cite{CC22}) using a scale-broken estimate
from the subelliptic theory of $D_{\xi }^{2}$ (see (4.7) of \cite{CC22} and
comments there).

To characterize $m(J,\theta )=0$ we need a trick, among others, inspired by
the idea of Schoen and Yau \cite{SY79} to show the torsion vanishes (see 
\cite[Lemma 4.7]{CC22}). To prove Corollary \ref{PMT'} above, we first blow
up the closed $M$ at a point $p$ by the Green's function of the $CR$
invariant sublaplacian $-(2+\frac{2}{n})\Delta _{b}+W$ to get an
asymptotically flat pseudohermitian manifold $N.$ Then we can apply Theorem %
\ref{PMT} to obtain the conclusion.

To solve the $CR$ Yamabe minimizer problem, we need a test function estimate
(see \cite[Theorem 5.1 in Section 5]{CC22}). The idea was rooted in an
argument used by Schoen in \cite{Sc84} for the Riemannian case. For the $CR$
case, it was first treated by Z. Li (\cite{Li90}) in an unpublished draft.
We reorganize his construction of a family of test functions $\phi _{\beta }$
and clarify the arguments at some points so that the $CR$-Sobolev quotient $%
Q(\phi _{\beta })$ is less than $\mathcal{Y}(S^{2n+1},\hat{J})$ minus a
positive multiple of the p-mass modulo the terms of higher decay rate (see 
\cite[(5.2)]{CC22}). From \cite[Theorem 5.1]{CC22} the result below follows
easily.

\begin{theorem}
\label{YMP} Suppose that $(M,\xi )$ is a closed (compact with no boundary),
contact and spin manifold of dimension 5. Assume that $J$ is a spherical CR
structure on $(M,\xi )$ with $\mathcal{Y(}M,J)$ $>$ $0$. Then the CR Yamabe
minimizer problem is solvable, i.e. we can find a solution to (\ref{1-1})
with minimum energy.
\end{theorem}

In \cite[Section 6]{CC22} we show that the connected sum of finitely many
(duplication allowed) 5-manifolds chosen arbitrarily from the set consisting
of $S^{5}/\mathbb{Z}_{p},$ $p$ an odd integer, $S^{4}\times S_{(a)}^{1},$ $%
a>1,$ and $\mathbb{RP}^{5}$ $\sharp $ $\mathbb{RP}^{5}$ is still a closed,
contact spin 5-manifold which admit a $spherical$ $CR$ structure with
positive $CR$ Yamabe constant (see \cite[Proposition 6.3]{CC22} or the
detailed proof in \cite{CC19}).

\end{document}